\documentclass[12pt,a4paper]{article}

\usepackage{amsmath}
\usepackage{amssymb}
\usepackage{bbold}
\usepackage{graphicx}
\usepackage{pstricks}

\allowdisplaybreaks

\renewcommand{\leq}{\leqslant}

\renewcommand{\geq}{\geqslant}

\makeatletter
\newcommand*{\bigcorr@macro}[2]{\sbox{0}{\mbox{$#1($}}\dimen0=\ht0
                \advance\dimen0 by \dp0
                \multiply\dimen0 by #2 \divide\dimen0 by 100}
\newcommand*{\bigcorr@big}[2]{\mbox{$#1\left#2\bigcorr@macro{#1}{85}\vrule
                   height \dimen0 depth 0pt width 0pt\right.\n@space$}}
\newcommand*{\bigcorr@Big}[2]{\mbox{$#1\left#2\bigcorr@macro{#1}{115}\vrule
                   height \dimen0 depth 0pt width 0pt\right.\n@space$}}
\newcommand*{\bigcorr@bigg}[2]{\mbox{$#1\left#2\bigcorr@macro{#1}{145}\vrule
                   height \dimen0 depth 0pt width 0pt\right.\n@space$}}
\newcommand*{\bigcorr@Bigg}[2]{\mbox{$#1\left#2\bigcorr@macro{#1}{175}\vrule
                   height \dimen0 depth 0pt width 0pt\right.\n@space$}}
\DeclareRobustCommand*{\big}[1]{{\mathpalette\bigcorr@big{#1}}}
\DeclareRobustCommand*{\Big}[1]{{\mathpalette\bigcorr@Big{#1}}}
\DeclareRobustCommand*{\bigg}[1]{{\mathpalette\bigcorr@bigg{#1}}}
\DeclareRobustCommand*{\Bigg}[1]{{\mathpalette\bigcorr@Bigg{#1}}}

\makeatother

\newtheorem{theorem}{Theorem}
\newtheorem{lemma}[theorem]{Lemma}

\let\ds\displaystyle
\newcommand{\cN}{\mathcal{N}}

\newcommand{\RR}{\mathbb{R}}
\newcommand{\ZZ}{\mathbb{Z}}
\newcommand{\argsup}{\mathop{\rm argsup}\limits}
\newcommand{\argmax}{\mathop{\rm argmax}\limits}
\newcommand*{\abs}[1]{\left|#1\right|}
\newcommand{\eb}{\mathbf{E}}
\newcommand{\pb}{\mathbf{P}}
\newcommand{\1}{\mathbb{1}}

\begin{document}

\title{On Limiting Likelihood Ratio Processes\\
       of some Change-Point Type Statistical Models}
\author{Sergue\"{\i} \textsc{Dachian}\\
Laboratoire de Math\'ematiques\\
Universit\'e Blaise Pascal\\
63177 Aubi\`ere CEDEX, France\\
Serguei.Dachian@math.univ-bpclermont.fr}

\date{}
\maketitle

\begin{abstract}
Different change-point type models encountered in statistical inference for
stochastic processes give rise to different limiting likelihood ratio
processes. In this paper we consider two such likelihood ratios. The first one
is an exponential functional of a two-sided Poisson process driven by some
parameter, while the second one is an exponential functional of a two-sided
Brownian motion. We establish that for sufficiently small values of the
parameter, the Poisson type likelihood ratio can be approximated by the
Brownian type one. As a consequence, several statistically interesting
quantities (such as limiting variances of different estimators) related to the
first likelihood ratio can also be approximated by those related to the second
one. Finally, we discuss the asymptotics of the large values of the parameter
and illustrate the results by numerical simulations.
\end{abstract}



\section{Introduction}

Different change-point type models encountered in statistical inference for
stochastic processes give rise to different limiting likelihood ratio
processes. In this paper we consider two of these processes. The first one is
the random process $Z_\rho$ on $\RR$ defined by
\begin{equation}
\label{proc1}
\ln Z_\rho(x)=\begin{cases}
\vphantom{\Big)}\rho\,\Pi_+(x)-x, &\text{if } x\geq 0,\\
\vphantom{\Big)}-\rho\,\Pi_-(-x)-x, &\text{if } x\leq 0,\\
\end{cases}
\end{equation}
where $\rho>0$, and $\Pi_+$ and $\Pi_-$ are two independent Poisson processes
on $\RR_+$ with intensities $1/(e^\rho-1)$ and $1/(1-e^{-\rho})$ respectively.
We also consider the random variables
\begin{equation}
\label{vars1}
\zeta_\rho=\frac{\int_{\RR}x\,Z_\rho(x)\;dx}{\int_{\RR}\,Z_\rho(x)\;dx}
\quad\text{and}\quad\xi_\rho=\argsup_{x\in\RR}Z_\rho(x)
\end{equation}
related to this process, as well as to their second moments
$B_\rho=\eb\zeta_\rho^2$ and $M_\rho=\eb\xi_\rho^2$.

The process $Z_\rho$ (up to a linear time change) arises in some non-regular,
namely change-point type, statistical models as the limiting likelihood ratio
process, and the variables $\zeta_\rho$ and $\xi_\rho$ (up to a multiplicative
constant) as the limiting distributions of the Bayesian estimators and of the
maximum likelihood estimator respectively. In particular, $B_\rho$ and
$M_\rho$ (up to the square of the above multiplicative constant) are the
limiting variances of these estimators, and the Bayesian estimators being
asymptotically efficient, the ratio $E_\rho=B_\rho/M_\rho$ is the asymptotic
efficiency of the maximum likelihood estimator in these models.

The main such model is the below detailed model of i.i.d.\ observations in the
situation when their density has a jump (is discontinuous). Probably the first
general result about this model goes back to Chernoff and
Rubin~\cite{CR}. Later, it was exhaustively studied by Ibragimov and
Khasminskii in~\cite[Chapter~5]{IKh} $\bigl($see also their previous
works~\cite{IKh2} and~\cite{IKh4}$\bigr)$.

\bigskip\bigskip\noindent
\textbf{Model 1.} Consider the problem of estimation of the location parameter
$\theta$ based on the observation $X^n=(X_1,\ldots,X_n)$ of the i.i.d.\ sample
from the density $f(x-\theta)$, where the known function $f$ is smooth enough
everywhere except at~$0$, and in $0$ we have
$$
0\ne\lim_{x\uparrow0}f(x)=a\ne b=\lim_{x\downarrow0}f(x)\ne0.
$$

Denote $\pb_{\theta}^n$ the distribution (corresponding to the parameter
$\theta$) of the observation $X^n$. As $n\to\infty$, the normalized likelihood
ratio process of this model defined by
$$
Z_n(u)=\frac{d\pb_{\theta+\frac un}^n}{d\pb_{\theta}^n}(X^n)=
\prod_{i=1}^{n}\frac{f\left(X_i-\theta-\frac un\right)}{f(X_i-\theta)}
$$
converges weakly in the space $\mathcal{D}_0(-\infty ,+\infty)$ (the Skorohod
space of functions on $\RR$ without discontinuities of the second kind and
vanishing at infinity) to the process $Z_{a,b}$ on $\RR$ defined by
$$
\ln Z_{a,b}(u)=\begin{cases}
\vphantom{\bigg)}\ln(\frac ab)\,\Pi_b(u)-(a-b)\,u, &\text{if } u\geq 0,\\
\vphantom{\bigg)}-\ln(\frac ab)\,\Pi_a(-u)-(a-b)\,u, &\text{if } u\leq 0,\\
\end{cases}
$$
where $\Pi_b$ and $\Pi_a$ are two independent Poisson processes on $\RR_+$
with intensities $b$ and $a$ respectively. The limiting distributions of the
Bayesian estimators and of the maximum likelihood estimator are given by
$$
\zeta_{a,b}=\frac{\int_{\RR}u\,Z_{a,b}(u)\;du}{\int_{\RR}\,Z_{a,b}(u)\;du}
\quad\text{and}\quad\xi_{a,b}=\argsup_{u\in\RR}Z_{a,b}(u)
$$
respectively. The convergence of moments also holds, and the Bayesian
estimators are asymptotically efficient. So, $\eb\zeta_{a,b}^2$ and
$\eb\xi_{a,b}^2$ are the limiting variances of these estimators, and
$\eb\zeta_{a,b}^2/\eb\xi_{a,b}^2$ is the asymptotic efficiency of the maximum
likelihood estimator.

Now let us note, that up to a linear time change, the process $Z_{a,b}$ is
nothing but the process $Z_\rho$ with $\rho=\abs{\ln(\frac ab)}$. Indeed, by
putting $u=\frac {x}{a-b}$ we get
\begin{align*}
\ln Z_{a,b}(u)&=\begin{cases}
\vphantom{\bigg)}\ln(\frac ab)\,\Pi_b(\frac {x}{a-b})-x, &\text{if }
\frac {x}{a-b}\geq 0,\\
\vphantom{\bigg)}-\ln(\frac ab)\,\Pi_a(-\frac {x}{a-b})-x, &\text{if }
\frac {x}{a-b}\leq 0,\\
\end{cases}\\
&=\ln Z_\rho(x)=\ln Z_\rho\bigl((a-b)\,u\bigr).
\end{align*}
So, we have
$$
\zeta_{a,b}=\frac{\zeta_\rho}{a-b}\quad\text{and}\quad
\xi_{a,b}=\frac{\xi_\rho}{a-b}\,,
$$
and hence
$$
\eb\zeta_{a,b}^2=\frac{B_\rho}{(a-b)^2}\;,\quad
\eb\xi_{a,b}^2=\frac{M_\rho}{(a-b)^2}\quad\text{and}\quad
\frac{\eb\zeta_{a,b}^2}{\eb\xi_{a,b}^2}=E_\rho\,.
$$

\bigskip\bigskip

Some other models where the process $Z_\rho$ arises occur in the statistical
inference for inhomogeneous Poisson processes, in the situation when their
intensity function has a jump (is discontinuous). In
Kutoyants~\cite[Chapter~5]{Kut} $\bigl($see also his previous
work~\cite{Kut2}$\bigr)$ one can find several examples, one of which is
detailed below.

\bigskip\bigskip\noindent
\textbf{Model 2.} Consider the problem of estimation of the location parameter
$\theta\in\left]\alpha,\beta\right[$, $0<\alpha<\beta<\tau$, based on the
observation $X^T$ on $[0,T]$ of the Poisson process with $\tau$-periodic
strictly positive intensity function $S(t+\theta)$, where the known function
$S$ is smooth enough everywhere except at points $t^*+\tau k$, $k\in\ZZ$, with
some $t^*\in\left[0,\tau\right]$, in which we have
$$
0\ne\lim_{t\uparrow t^*}S(t)=S_-\ne S_+=\lim_{t\downarrow t^*}S(t)\ne0.
$$

Denote $\pb_{\theta}^T$ the distribution (corresponding to the parameter
$\theta$) of the observation $X^T$. As $T\to\infty$, the normalized likelihood
ratio process of this model defined by
$$
Z_T(u)=\frac{d\pb_{\theta+\frac uT}^T}{d\pb_{\theta}^T}(X^T)=\exp
\biggl\{\int_{0}^{T}\!\!\ln\frac{S_{\theta+\frac uT}(t)}{S_{\theta}(t)}\,dX(t)
-\int_{0}^{T}\!\!\left[S_{\theta+\frac uT}(t)-S_\theta(t)\right]\!dt\biggr\}
$$
converges weakly in the space $\mathcal{D}_0(-\infty ,+\infty)$ to the process
$Z_{\tau,S_-,S_+}$ on $\RR$ defined by
$$
\ln Z_{\tau,S_-,S_+}=\begin{cases}
\vphantom{\bigg)}\ln\bigl(\frac{S_+}{S_-}\bigr)\,\Pi_{S_-}\bigl(\frac
u\tau\bigr)-(S_+-S_-)\,\frac u\tau\,, &\text{if } u\geq 0,\\
\vphantom{\bigg)}-\ln\bigl(\frac{S_+}{S_-}\bigr)\,\Pi_{S_+}\bigl(-\frac
u\tau\bigr)-(S_+-S_-)\,\frac u\tau\,, &\text{if } u\leq 0,\\
\end{cases}
$$
where $\Pi_{S_-}$ and $\Pi_{S_+}$ are two independent Poisson processes on
$\RR_+$ with intensities $S_-$ and $S_+$ respectively. The limiting
distributions of the Bayesian estimators and of the maximum likelihood
estimator are given by
$$
\zeta_{\tau,S_-,S_+}=\frac{\int_{\RR}u\,Z_{\tau,S_-,S_+}(u)\;du}
{\int_{\RR}\,Z_{\tau,S_-,S_+}(u)\;du}\quad\text{and}\quad
\xi_{\tau,S_-,S_+}=\argsup_{u\in\RR}Z_{\tau,S_-,S_+}(u)
$$
respectively. The convergence of moments also holds, and the Bayesian
estimators are asymptotically efficient. So, $\eb\zeta_{\tau,S_-,S_+}^2$ and
$\eb\xi_{\tau,S_-,S_+}^2$ are the limiting variances of these estimators, and
$\eb\zeta_{\tau,S_-,S_+}^2/\eb\xi_{\tau,S_-,S_+}^2$ is the asymptotic
efficiency of the maximum likelihood estimator.

Now let us note, that up to a linear time change, the process
$Z_{\tau,S_-,S_+}$ is nothing but the process $Z_\rho$ with
$\rho=\bigl|\ln\bigl(\frac{S_+}{S_-}\bigr)\bigr|$. Indeed, by putting
$u=\frac{\tau x}{S_+-S_-}$ we get
$$
Z_{\tau,S_-,S_+}(u)=Z_\rho(x)= Z_\rho\left(\frac{S_+-S_-}{\tau}\,u\right).
$$
So, we have
$$
\zeta_{\tau,S_-,S_+}=\frac{\tau\,\zeta_\rho}{S_+-S_-}\quad\text{and}\quad
\zeta_{\tau,S_-,S_+}= \frac{\tau\,\xi_\rho}{S_+-S_-}\,,
$$
and hence
$$
\eb\zeta_{\tau,S_-,S_+}^2=\frac{\tau^2\,B_\rho}{(S_+-S_-)^2}\;,\quad
\eb\xi_{\tau,S_-,S_+}^2=\frac{\tau^2\,M_\rho}{(S_+-S_-)^2}\quad\text{and}\quad
\frac{\eb\zeta_{\tau,S_-,S_+}^2}{\eb\xi_{\tau,S_-,S_+}^2}=E_\rho\,.
$$

\bigskip\bigskip

The second limiting likelihood ratio process considered in this paper is the
random process
\begin{equation}
\label{proc2}
Z_0(x)=\exp\left\{W(x)-\frac12\abs x\right\},\quad x\in\RR,
\end{equation}
where $W$ is a standard two-sided Brownian motion. In this case, the limiting
distributions of the Bayesian estimators and of the maximum likelihood
estimator (up to a multiplicative constant) are given by
\begin{equation}
\label{vars2}
\zeta_0=\frac{\int_{\RR}x\,Z_0(x)\;dx}{\int_{\RR}\,Z_0(x)\;dx}
\quad\text{and}\quad\xi_0=\argsup_{x\in\RR}Z_0(x)
\end{equation}
respectively, and the limiting variances of these estimators (up to the square
of the above multiplicative constant) are $B_0=\eb\zeta_0^2$ and
$M_0=\eb\xi_0^2$.

The models where the process $Z_0$ arises occur in various fields of
statistical inference for stochastic processes. A well-known example is the
below detailed model of a discontinuous signal in a white Gaussian noise
exhaustively studied by Ibragimov and Khasminskii in~\cite[Chapter~7.2]{IKh}
$\bigl($see also their previous work~\cite{IKh5}$\bigr)$, but one can also
cite change-point type models of dynamical systems with small noise
$\bigl($see Kutoyants~\cite{Kut2} and~\cite[Chapter~5]{Kut3}$\bigr)$, those of
ergodic diffusion processes $\bigl($see
Kutoyants~\cite[Chapter~3]{Kut4}$\bigr)$, a change-point type model of delay
equations $\bigl($see K\"uchler and Kutoyants~\cite{KK}$\bigr)$, an
i.i.d.\ change-point type model $\bigl($see Deshayes and
Picard~\cite{DP}$\bigr)$, a model of a discontinuous periodic signal in a time
inhomogeneous diffusion $\bigl($see H\"opfner and Kutoyants~\cite{HK}$\bigr)$,
and so on.

\bigskip\bigskip\noindent
\textbf{Model 3.} Consider the problem of estimation of the location parameter
$\theta\in\left]\alpha,\beta\right[$, $0<\alpha<\beta<1$, based on the
observation $X^\varepsilon$ on $\left[0,1\right]$ of the random process
satisfying the stochastic differential equation
$$
dX^\varepsilon(t)=\frac1\varepsilon\,S(t-\theta)\,dt+dW(t),
$$
where $W$ is a standard Brownian motion, and $S$ is a known function having a
bounded derivative on $\left]-1,0\right[\cup\left]0,1\right[$ and satisfying
$$
\lim_{t\uparrow0}S(t)-\lim_{t\downarrow0}S(t)=r\ne0.
$$

Denote $\pb_{\theta}^\varepsilon$ the distribution (corresponding to the
parameter $\theta$) of the observation $X^\varepsilon$. As $\varepsilon\to0$,
the normalized likelihood ratio process of this model defined by
\begin{align*}
Z_\varepsilon(u)&=\frac{d\pb_{\theta+\varepsilon^2u}^\varepsilon}
{d\pb_{\theta}^\varepsilon}(X^\varepsilon)\\
&=\exp\biggl\{\frac1\varepsilon\int_{0}^{1}
\bigl[S(t-\theta-\varepsilon^2u)-S(t-\theta)\bigr]\,dW(t)\\
&\phantom{=\exp\biggl\{}-\frac1{2\,\varepsilon^2}\int_{0}^{1}
\bigl[S(t-\theta-\varepsilon^2u)-S(t-\theta)\bigr]^2\,dt \biggr\}
\end{align*}
converges weakly in the space $\mathcal{C}_0(-\infty ,+\infty)$ (the space of
continuous functions vanishing at infinity equipped with the supremum norm) to
the process $Z_0(r^2u)$, $u\in\RR$. The limiting distributions of the Bayesian
estimators and of the maximum likelihood estimator are $r^{-2}\zeta_0$ and
$r^{-2}\xi_0$ respectively. The convergence of moments also holds, and the
Bayesian estimators are asymptotically efficient. So, $r^{-4}B_0$ and
$r^{-4}M_0$ are the limiting variances of these estimators, and $E_0$ is the
asymptotic efficiency of the maximum likelihood estimator.

\bigskip\bigskip

Let us also note that Terent'yev in~\cite{Ter} determined explicitly the
distribution of $\xi_0$ and calculated the constant $M_0=26$. These results
were taken up by Ibragimov and Khasminskii in~\cite[Chapter~7.3]{IKh}, where
by means of numerical simulation they equally showed that $B_0=19.5\pm0.5$,
and so $E_0=0.73\pm0.03$. Later in~\cite{Gol}, Golubev expressed $B_0$ in
terms of the second derivative (with respect to a parameter) of an improper
integral of a composite function of modified Hankel and Bessel
functions. Finally in~\cite{RS}, Rubin and Song obtained the exact values
$B_0=16\,\zeta(3)$ and $E_0=8\,\zeta(3)/13$, where~$\zeta$ is Riemann's zeta
function defined by
$$
\zeta(s)=\sum_{n=1}^{\infty}\frac1{n^s}\;.
$$

The random variables $\zeta_\rho$ and $\xi_\rho$ and the quantities $B_\rho$,
$M_\rho$ and $E_\rho$, $\rho>0$, are much less studied. One can cite
Pflug~\cite{Pfl} for some results about the distribution of the random
variables
$$
\argsup_{x\in\RR_+}Z_\rho(x)\quad\text{and}\quad\argsup_{x\in\RR_-}Z_\rho(x)
$$
related to $\xi_\rho$.

In this paper we establish that the limiting likelihood ratio processes
$Z_\rho$ and $Z_0$ are related. More precisely, we show that as $\rho\to 0$,
the process $Z_\rho(y/\rho)$, $y\in\RR$, converges weakly in the space
$\mathcal{D}_0(-\infty ,+\infty)$ to the process~$Z_0$. So, the random
variables $\rho\,\zeta_\rho$ and $\rho\,\xi_\rho$ converge weakly to the
random variables $\zeta_0$ and $\xi_0$ respectively. We show equally that the
convergence of moments of these random variables holds, that is, $\rho^2
B_\rho \to 16\,\zeta(3)$, $\rho^2 M_\rho \to 26$ and $E_\rho \to
8\,\zeta(3)/13$.

These are the main results of the present paper, and they are presented in
Section~\ref{MR}, where we also briefly discuss the second possible
asymptotics $\rho\to+\infty$. The necessary lemmas are proved in
Section~\ref{PL}. Finally, some numerical simulations of the quantities
$B_\rho$, $M_\rho$ and $E_\rho$ for $\rho\in\left]0,\infty\right[$ are
presented in Section~\ref{NS}.

\section{Main results}
\label{MR}

Consider the process $X_\rho(y)=Z_\rho(y/\rho)$, $y\in\RR$, where $\rho>0$ and
$Z_\rho$ is defined by~\eqref{proc1}. Note that
$$
\frac{\int_{\RR}y\,X_\rho(y)\;dy}{\int_{\RR}\,X_\rho(y)\;dy}=\rho \,\zeta_\rho
\quad\text{and}\quad\argsup_{y\in\RR}X_\rho(y)=\rho\,\xi_\rho\,,
$$
where the random variables $\zeta_\rho$ and $\xi_\rho$ are defined
by~\eqref{vars1}. Remind also the process $Z_0$ on $\RR$ defined
by~\eqref{proc2} and the random variables $\zeta_0$ and $\xi_0$ defined
by~\eqref{vars2}. Recall finally the quantities $B_\rho=\eb\zeta_\rho^2$,
$M_\rho=\eb\xi_\rho^2$, $E_\rho=B_\rho/M_\rho$,
$B_0=\eb\zeta_0^2=16\,\zeta(3)$, $M_0=\eb\xi_0^2=26$ and
$E_0=B_0/M_0=8\,\zeta(3)/13$. Now we can state the main result of the present
paper.

\begin{theorem}
\label{T1}
The process $X_\rho$ converges weakly in the space $\mathcal{D}_0(-\infty
,+\infty)$ to the process $Z_0$ as $\rho \to 0$. In particular, the random
variables $\rho\,\zeta_\rho$ and $\rho\,\xi_\rho$ converge weakly to the
random variables $\zeta_0$ and $\xi_0$ respectively. Moreover, for any $k>0$
we have
$$
\rho^k\,\eb\zeta_\rho^k \to \eb\zeta_0^k\quad\text{and}\quad
\rho^k\,\eb\xi_\rho^k \to \eb\xi_0^k,
$$
and in particular $\rho^2 B_\rho \to 16\,\zeta(3)$, $\rho^2 M_\rho \to 26$ and
$E_\rho \to 8\,\zeta(3)/13$.
\end{theorem}

The results concerning the random variable $\zeta_\rho$ are direct consequence
of Ibragimov and Khasminskii~\cite[Theorem~1.10.2]{IKh} and the following
three lemmas.

\begin{lemma}
\label{L1}
The finite-dimensional distributions of the process $X_\rho$ converge to those
of $Z_0$ as $\rho \to 0$.
\end{lemma}

\begin{lemma}
\label{L2}
For all $\rho>0$ and all $y_1,y_2\in\RR$ we have
$$
\eb\left|X_\rho^{1/2}(y_1)-X_\rho^{1/2}(y_2)\right|^2\leq
\frac14\abs{y_1-y_2}.
$$
\end{lemma}

\begin{lemma}
\label{L3}
For any $c\in\left]\,0\,{,}\;1/8\,\right[$ we have
$$
\eb X_\rho^{1/2}(y)\leq\exp\bigl(-c\abs y\bigr)
$$
for all sufficiently small $\rho$ and all $y\in\RR$.
\end{lemma}

Note that these lemmas are not sufficient to establish the weak convergence of
the process $X_\rho$ in the space $\mathcal{D}_0(-\infty ,+\infty)$ and the
results concerning the random variable $\xi_\rho$. However, the increments of
the process $\ln X_\rho$ being independent, the convergence of its
restrictions (and hence of those of $X_\rho$) on finite intervals
$[A,B]\subset\RR$ $\bigl($that is, convergence in the Skorohod space
$\mathcal{D}[A,B]$ of functions on $[A,B]$ without discontinuities of the
second kind$\bigr)$ follows from Gihman and Skorohod~\cite[Theorem~6.5.5]{GS},
Lemma~\ref{L1} and the following lemma.

\begin{lemma}
\label{L4}
For any $\varepsilon>0$ we have
$$
\lim_{h\to 0}\ \lim_{\rho\to 0}\ \sup_{\abs{y_1-y_2} < h}
\pb\Bigl\{\bigl|\ln X_\rho(y_1)-\ln X_\rho(y_2)\bigr|>\varepsilon\Bigr\}=0.
$$
\end{lemma}

Now, Theorem~\ref{T1} follows from the following estimate on the tails of the
process $X_\rho$ by standard argument.

\begin{lemma}
\label{L5}
For any $b\in\left]\,0\,{,}\;3/40\,\right[$ we have
$$
\pb\biggl\{\sup_{\abs y>A} X_\rho(y) > e^{-bA}\biggr\} \leq 2\,e^{-bA}
$$
for all sufficiently small $\rho$ and all $A>0$.
\end{lemma}

All the above lemmas will be proved in the next section, but before let us
discuss the second possible asymptotics $\rho\to+\infty$. One can show that in
this case, the process $Z_\rho$ converges weakly in the space
$\mathcal{D}_0(-\infty ,+\infty)$ to the process
$Z_\infty(u)=e^{-u}\,\1_{\{u>\eta\}}$, $u\in\RR$, where $\eta$ is a negative
exponential random variable with $\pb\{\eta<t\}=e^t$, $t\leq 0$. So, the
random variables $\zeta_\rho$ and $\xi_\rho$ converge weakly to the random
variables
$$
\zeta_\infty=
\frac{\int_{\RR}u\,Z_\infty(u)\;du}{\int_{\RR}\,Z_\infty(u)\;du}=\eta+1
\quad\text{and}\quad \xi_\infty=\argsup_{u\in\RR}Z_\infty(u)=\eta
$$
respectively. One can equally show that, moreover, for any $k>0$ we have
$$
\eb\zeta_\rho^k\to\eb\zeta_\infty^k \quad\text{and}\quad
\eb\xi_\rho^k\to\eb\xi_\infty^k,
$$
and in particular, denoting $B_\infty\!=\eb\zeta_\infty^2$,
$M_\infty\!=\eb\xi_\infty^2$ and $E_\infty\!=B_\infty/M_\infty$, we finally
have $B_\rho \to B_\infty\!=\eb(\eta+1)^2=1$, $M_\rho \to
M_\infty\!=\eb\eta^2=2$ and $E_\rho \to E_\infty\!=1/2$.

Let us note that these convergences are natural, since the process $Z_\infty$
can be considered as a particular case of the process $Z_\rho$ with
$\rho=+\infty$ if one admits the convention $+\infty\cdot0=0$.

Note also that the process $Z_\infty$ (up to a linear time change) is the
limiting likelihood ratio process of Model~1 (Model~2) in the situation when
$a\cdot b=0$ ($S_-\cdot S_+$=0). In this case, the variables
$\zeta_\infty=\eta+1$ and $\xi_\infty=\eta$ (up to a multiplicative constant)
are the limiting distributions of the Bayesian estimators and of the maximum
likelihood estimator respectively. In particular, $B_\infty=1$ and
$M_\infty=2$ (up to the square of the above multiplicative constant) are the
limiting variances of these estimators, and the Bayesian estimators being
asymptotically efficient, $E_\infty=1/2$ is the asymptotic efficiency of the
maximum likelihood estimator.

\section{Proofs of the lemmas}
\label{PL}

First we prove Lemma~\ref{L1}. Note that the restrictions of the process $\ln
X_\rho$ (as well as those of the process $\ln Z_0$) on $\RR_+$ and on $\RR_-$
are mutually independent processes with stationary and independent increments.
So, to obtain the convergence of all the finite-dimensional distributions, it
is sufficient to show the convergence of one-dimensional distributions only,
that is,
$$
\ln X_\rho(y)\Rightarrow\ln Z_0(y)=W(y)-\frac{\abs{y}}{2}=\cN
\Bigl(-\frac{\abs{y}}{2}\,,\abs{y}\Bigr)
$$
for all $y\in\RR$. Here and in the sequel ``$\Rightarrow$'' denotes the weak
convergence of the random variables, and $\cN(m,V)$ denotes a ``generic''
random variable distributed according to the normal law with mean $m$ and
variance $V$.

Let $y>0$.  Then, noting that $\ds\Pi_+\Bigl(\frac{y}{\rho}\Bigr)$ is a
Poisson random variable of parameter
$\ds\lambda=\frac{y}{\rho\,(e^\rho-1)}\to\infty$, we have
\begin{align*}
\ln X_\rho(y)&=\rho\,\Pi_+\Bigl(\frac{y}{\rho}\Bigr)-\frac{y}{\rho}
=\rho\,\sqrt{\frac{y}{\rho\,(e^\rho-1)}}\ \frac
{\Pi_+\bigl(\frac{y}{\rho}\bigr)-\lambda}{\sqrt{\lambda}}
+\frac{y}{e^\rho-1}-\frac{y}{\rho}\\
&=\sqrt{y}\,\sqrt{\frac{\rho}{e^\rho-1}}\ \frac
{\Pi_+\bigl(\frac{y}{\rho}\bigr)-\lambda}{\sqrt{\lambda}}
-y\,\frac{e^\rho-1-\rho}{\rho\,(e^\rho-1)}\Rightarrow
\cN\left(-\frac{y}{2}\,,y\right),
\end{align*}
since
$$
\frac{\rho}{e^\rho-1}=\frac{\rho}{\rho+o(\rho)}\to1,\qquad
\frac{e^\rho-1-\rho}{\rho\,(e^\rho-1)}=
\frac{\rho^2/2+o(\rho^2)}{\rho\,\bigr(\rho+o(\rho)\bigr)}\to\frac{1}{2}
$$
and
$$
\frac{\Pi_+\bigl(\frac{y}{\rho}\bigr)-\lambda}{\sqrt{\lambda}}\Rightarrow
\cN(0,1).
$$

Similarly, for $y<0$ we have
\begin{align*}
\ln X_\rho(y)&=-\rho\,\Pi_-\Bigl(\frac{-y}{\rho}\Bigr)-\frac{y}{\rho}
=\rho\,\sqrt{\frac{-y}{\rho\,(1-e^{-\rho})}}\ \frac
{\lambda'-\Pi_-\bigl(\frac{-y}{\rho}\bigr)}{\sqrt{\lambda'}}
-\frac{-y}{1-e^{-\rho}}-\frac{y}{\rho}\\
&=\sqrt{-y}\,\sqrt{\frac{\rho}{1-e^{-\rho}}}\ \frac
{\lambda'-\Pi_-\bigl(\frac{-y}{\rho}\bigr)}{\sqrt{\lambda'}}
+y\,\frac{e^{-\rho}-1+\rho}{\rho\,(1-e^{-\rho})}\Rightarrow
\cN\left(\frac{y}{2}\,,-y\right),
\end{align*}
and so, Lemma~\ref{L1} is proved.

\bigskip\bigskip

Now we turn to the proof of Lemma~\ref{L3} (we will prove Lemma~\ref{L2}
just after). For $y>0$ we can write
$$
\eb X_\rho^{1/2}(y)=\eb\exp\left(\frac\rho2\,\Pi_+\Bigl(\frac
y\rho\Bigr)-\frac y{2\rho}\right)=\exp\left(-\frac y{2\rho}\right)\,
\eb\exp\left(\frac\rho2\,\Pi_+\Bigl(\frac y\rho\Bigr)\right).
$$
Note that $\ds\Pi_+\Bigl(\frac{y}{\rho}\Bigr)$ is a Poisson random variable of
parameter $\ds\lambda=\frac{y}{\rho\,(e^\rho-1)}$ with moment generating
function $M(t)=\exp\bigl(\lambda\,(e^t-1)\bigr)$. So, we get
\begin{align*}
\eb X_\rho^{1/2}(y)&=\exp\left(-\frac y{2\rho}\right)\,\exp\left(\frac
y{\rho\,(e^\rho-1)}\,(e^{\rho/2}-1)\right)\\
&=\exp\left(-\frac y{2\rho}+\frac y{\rho\,(e^{\rho/2}+1)}\right)
=\exp\left(-y\,\frac{e^{\rho/2}-1}{2\rho\,(e^{\rho/2}+1)}\right)\\
&=\exp\left(-y\,\frac{e^{\rho/4}-e^{-\rho/4}}
{2\rho\,(e^{\rho/4}+e^{-\rho/4})}\right)=\exp
\left(-y\,\frac{\tanh(\rho/4)}{2\rho}\right).
\end{align*}

For $y<0$ we obtain similarly
\begin{align*}
\eb X_\rho^{1/2}(y)&=\eb\exp\left(-\frac\rho2\,\Pi_-\Bigl(\frac
{-y}\rho\Bigr)-\frac y{2\rho}\right)\\
&=\exp\left(-\frac y{2\rho}\right)\,\exp
\left(\frac{-y}{\rho\,(1-e^{-\rho})}\,(e^{-\rho/2}-1)\right)\\
&=\exp\left(-\frac y{2\rho}+\frac y{\rho\,(1+e^{-\rho/2})}\right)=\exp
\left(y\,\frac{1-e^{-\rho/2}}{2\rho\,(1+e^{-\rho/2})}\right)\\
&=\exp\left(y\,\frac{\tanh(\rho/4)}{2\rho}\right).
\end{align*}

Thus, for all $y\in\RR$ we have
\begin{equation}
\label{Eroot}
\eb X_\rho^{1/2}(y)=\exp\left(-\abs y\frac{\tanh(\rho/4)}{2\rho}\right),
\end{equation}
and since
$$
\frac{\tanh(\rho/4)}{2\rho}=\frac{\rho/4+o(\rho)}{2\rho}\to\frac18
$$
as $\rho\to0$, for any $c\in\left]\,0\,{,}\;1/8\,\right[$ we have $\eb
X_\rho^{1/2}(y)\leq\exp\bigl(-c\abs{y}\bigr)$ for all sufficiently small
$\rho$ and all $y\in\RR$. Lemma~\ref{L3} is proved.

\bigskip\bigskip

Further we verify Lemma~\ref{L2}. We first consider the case $y_1,y_2\in\RR_+$
(say $y_1\geq y_2$). Using~\eqref{Eroot} and taking into account the
stationarity and the independence of the increments of the process $\ln
X_\rho$ on $\RR_+$, we can write
\begin{align*}%
\eb\left|X_\rho^{1/2}(y_1)-X_\rho^{1/2}(y_2)\right|^2&=\eb X_\rho(y_1)+\eb
X_\rho(y_2)-2\, \eb X_\rho^{1/2}(y_1)X_\rho^{1/2}(y_2)\\
&=2-2\,\eb X_\rho(y_2)\,\eb\frac{X_\rho^{1/2}(y_1)}{X_\rho^{1/2}(y_2)}\\
&=2-2\,\eb X_\rho^{1/2}(y_1-y_2)\\
&=2-2\exp\left(-\abs{y_1-y_2}\frac{\tanh(\rho/4)}{2\rho}\right)\\
&\leq\abs{y_1-y_2}\frac{\tanh(\rho/4)}{\rho}\leq\frac14\abs{y_1-y_2}.
\end{align*}

The case $y_1,y_2\in\RR_-$ can be treated similarly.

Finally, if $y_1y_2\leq 0$ (say $y_2\leq 0\leq y_1$), we have
\begin{align*}%
\eb\left|X_\rho^{1/2}(y_1)-X_\rho^{1/2}(y_2)\right|^2&=2-2\,\eb
X_\rho^{1/2}(y_1)\, \eb X_\rho^{1/2}(y_2)\\
&=2-2\exp\left(-\abs{y_1}\frac{\tanh(\rho/4)}{2\rho}-
\abs{y_2}\frac{\tanh(\rho/4)}{2\rho}\right)\\
&=2-2\exp\left(-\abs{y_1-y_2}\frac{\tanh(\rho/4)}{2\rho}\right)\\
&\leq\frac14\abs{y_1-y_2},
\end{align*}
and so, Lemma~\ref{L2} is proved.

\bigskip\bigskip

Now let us check Lemma~\ref{L4}. First let $y_1,y_2\in\RR_+$ (say $y_1\geq
y_2$) such that $\Delta=\abs{y_1-y_2}<h$. Then
\begin{align*}
\pb\Bigl\{\bigl|\ln X_\rho(y_1)-\ln X_\rho(y_2)\bigr|>\varepsilon\Bigr\}
&\leq\frac1{\varepsilon^2}\,\eb\bigl|\ln X_\rho(y_1)-\ln X_\rho(y_2)\bigr|^2\\
&=\frac1{\varepsilon^2}\,\eb\bigl|\ln X_\rho(\Delta)\bigr|^2\\
&=\frac1{\varepsilon^2}\,\eb\left|\rho\,\Pi_+
\Bigl(\frac{\Delta}{\rho}\Bigr)-\frac{\Delta}{\rho}\right|^2\\
&=\frac1{\varepsilon^2}\left(\rho^2(\lambda+\lambda^2)
+\frac{\Delta^2}{\rho^2}-2\lambda\Delta\right)\\
&=\frac1{\varepsilon^2}\left(\beta(\rho)\,\Delta
+\gamma(\rho)\,\Delta^2\right)\\
&<\frac1{\varepsilon^2}\left(\beta(\rho)\,h+\gamma(\rho)\,h^2\right),
\end{align*}
where $\lambda=\frac{\Delta}{\rho\,(e^\rho-1)}$ is the parameter of the
Poisson random variable $\Pi_+\bigl(\!\frac{\Delta}{\rho}\!\bigr)$,
\begin{align*}
\beta(\rho)&=\frac{\rho}{(e^\rho-1)}=\frac{\rho}{\rho+o(\rho)}\to1\\
\intertext{and}
\gamma(\rho)&=\frac{1}{(e^\rho-1)^2}+\frac{1}{\rho^2}-\frac{2}{\rho\,(e^\rho-1)}
=\left(\frac{1}{\rho}-\frac{1}{e^\rho-1}\right)^2\\
&=\biggl(\frac{e^\rho-1-\rho}{\rho\,(e^\rho-1)}\biggr)^2=\biggl(
\frac{\rho^2/2+o(\rho^2)}{\rho\,\bigl(\rho+o(\rho)\bigr)}\biggr)^2\to\frac14
\end{align*}
as $\rho\to 0$. So, we have
\begin{align*}
\lim_{\rho\to 0}\ \sup_{\abs{y_1-y_2} < h}\pb\Bigl\{\bigl|\ln X_\rho(y_1)-\ln
X_\rho(y_2)\bigr|>\varepsilon\Bigr\} &\leq\lim_{\rho\to 0}
\frac1{\varepsilon^2}\left(\beta(\rho)\,h+\gamma(\rho)\,h^2\right)\\
&=\frac1{\varepsilon^2}\left(h+\frac{h^2}4\right),
\end{align*}
and hence
$$
\lim_{h\to 0}\ \lim_{\rho\to 0}\ \sup_{\abs{y_1-y_2} < h}\pb\Bigl\{\bigl|\ln
X_\rho(y_1)-\ln X_\rho(y_2)\bigr|>\varepsilon\Bigr\}=0,
$$
where the supremum is taken only over $y_1,y_2\in\RR_+$.

For $y_1,y_2\in\RR_-$ such that $\Delta=\abs{y_1-y_2}<h$ one can obtain
similarly
\begin{align*}
\pb\Bigl\{\bigl|\ln X_\rho(y_1)-\ln X_\rho(y_2)\bigr|>\varepsilon\Bigr\}
&\leq\frac1{\varepsilon^2}\,\eb\bigl|\ln X_\rho(y_1)-\ln X_\rho(y_2)\bigr|^2\\
&=\frac1{\varepsilon^2}\left(\beta'(\rho)\,\Delta
+\gamma'(\rho)\,\Delta^2\right)\\
&<\frac1{\varepsilon^2}\left(\beta'(\rho)\,h+\gamma'(\rho)\,h^2\right),
\end{align*}
where
\begin{align*}
\beta'(\rho)&=\frac{\rho}{(1-e^{-\rho})}=\frac{\rho}{\rho+o(\rho)}\to1\\
\intertext{and}
\gamma'(\rho)&=\biggl(\frac{e^{-\rho}-1+\rho}{\rho\,(1-e^{\rho})}\biggr)^2=
\biggl(\frac{\rho^2/2+o(\rho^2)}{\rho\,\bigl(\rho+o(\rho)\bigr)}\biggr)^2\to
\frac14
\end{align*}
as $\rho\to 0$, which will yield the same conclusion as above, but with the
supremum taken over $y_1,y_2\in\RR_-$.

Finally, for $y_1y_2\leq 0$ (say $y_2\leq 0\leq y_1$) such that
$\abs{y_1-y_2}<h$, using the elementary inequality $(a-b)^2\leq 2(a^2+b^2)$ we
get
\begin{align*}
\pb\Bigl\{\bigl|\ln X_\rho(y_1)-\ln X_\rho(y_2)\bigr|>\varepsilon\Bigr\}
&\leq\frac1{\varepsilon^2}\,\eb\bigl|\ln X_\rho(y_1)-\ln X_\rho(y_2)\bigr|^2\\
&\leq\frac2{\varepsilon^2}\left(\eb\bigl|\ln X_\rho(y_1)\bigr|^2+
\eb\bigl|\ln X_\rho(y_2)\bigr|^2\right)\\
&=\frac2{\varepsilon^2}\left(\beta(\rho)y_1\!+\!\gamma(\rho)y_1^2
\!+\!\beta'(\rho)\!\abs{y_2}\!+\!\gamma'(\rho)\!\abs{y_2}^2\right)\\
&\leq\frac2{\varepsilon^2}\Bigl(\bigl(\beta(\rho)+\beta'(\rho)\bigr)\,h
+\bigl(\gamma(\rho)+\gamma'(\rho)\bigr)\,h^2\Bigr),
\end{align*}
which again will yield the desired conclusion. Lemma~\ref{L4} is proved.

\bigskip\bigskip

It remains to verify Lemma~\ref{L5}. Clearly,
$$
\pb\biggl\{\sup_{\abs y>A} X_\rho(y) > e^{-bA}\biggr\}\leq
\pb\biggl\{\sup_{y>A} X_\rho(y) > e^{-bA}\biggr\}+
\pb\biggl\{\sup_{y<-A} X_\rho(y) > e^{-bA}\biggr\}.
$$
In order to estimate the first term, we need two auxiliary results.

\begin{lemma}
\label{L6}
For any $c\in\left]\,0\,{,}\;3/32\,\right[$ we have
$$
\eb X_\rho^{1/4}(y)\leq\exp\bigl(-c\abs y\bigr)
$$
for all sufficiently small $\rho$ and all $y\in\RR$.
\end{lemma}

\begin{lemma}
\label{L7}
For all $\rho>0$ the random variable
$$
\eta_\rho=\sup_{t\in\RR_+}\bigl(\Pi_\lambda(t)-t\bigr),
$$
where $\Pi_\lambda$ is a Poisson process on $\RR_+$ with intensity
$\lambda=\rho/(e^\rho-1)\in\left]0,1\right[$, verifies
$$
\eb\exp\left(\frac\rho4\,\eta_\rho\right)\leq 2.
$$
\end{lemma}

The first result can be easily obtained following the proof of Lemma~\ref{L3},
so we prove the second one only. For this, let us remind that according to
Shorack and Wellner~\cite[Proposition~1 on page~392]{ShW} $\bigl($see also
Pyke~\cite{Pyke}$\bigr)$, the distribution function
$F_\rho(x)=\pb\{\eta_\rho<x\}$ of $\eta_\rho$ is given by
$$
1-F_\rho(x)=\pb\{\eta_\rho\geq x\}=(1-\lambda)\,e^{\lambda x}
\sum_{n>x}\frac{(n-x)^n}{n!}\,\bigl(\lambda\,e^{-\lambda}\bigr)^n
$$
for $x>0$, and is zero for $x\leq 0$. Hence, for $x>0$ we have
\begin{align*}
1-F_\rho(x)&\leq(1-\lambda)\,e^{\lambda x}\,\sum_{n>x}
\frac{(n-x)^n}{\sqrt{2\pi n}\,n^n\,e^{-n}}\,
\bigl(\lambda\,e^{-\lambda}\bigr)^n\\
&=\frac{1-\lambda}{\sqrt{2\pi}}\,e^{\lambda x}\,\sum_{n>x}\frac1{\sqrt{n}}
\left(1-\frac{x}{n}\right)^n \bigl(\lambda\,e^{1-\lambda}\bigr)^n\\
&\leq\frac{1-\lambda}{\sqrt{2\pi}}\,e^{\lambda x}\,\sum_{n>x}
e^{-x}\frac{\bigl(\lambda\,e^{1-\lambda}\bigr)^n}{\sqrt{n}}\\
&\leq\frac{1-\lambda}{\sqrt{2\pi}}\,e^{(\lambda-1)x}\,
\bigl(\lambda\,e^{1-\lambda}\bigr)^x\sum_{n>x}
\frac{\bigl(\lambda\,e^{1-\lambda}\bigr)^{n-x}}{\sqrt{n-x}}\\
&=\frac{1-\lambda}{\sqrt{2\pi}}\,\lambda^x
\sum_{k>0}\frac{\bigl(\lambda\,e^{1-\lambda}\bigr)^k}{\sqrt{k}}\leq
\frac{1-\lambda}{\sqrt{2\pi}}\,\lambda^x
\int_{\RR_+}\frac{\bigl(\lambda\,e^{1-\lambda}\bigr)^t}{\sqrt{t}}\;dt\\
&=\frac{1-\lambda}{\sqrt{2\pi}}\,\lambda^x
\,\frac{\Gamma(1/2)}{\sqrt{-\ln\bigl(\lambda\,e^{1-\lambda}\bigr)}}
=\frac{1-\lambda}{\sqrt{-2\ln\bigl(\lambda\,e^{1-\lambda}\bigr)}}
\left(\frac{\rho}{e^\rho-1}\right)^x\\
&\leq\left(\frac{\rho\,e^{-\rho/2}}{e^{\rho/2}-e^{-\rho/2}}\right)^x
=\left(\frac{\rho\,e^{-\rho/2}}{2\sinh(\rho/2)}\right)^x\leq e^{-\rho x/2},
\end{align*}
where we used Stirling inequality and the inequality
$1-\lambda\leq\sqrt{-2\ln\bigl(\lambda\,e^{1-\lambda}\bigr)}$, which is easily
reduced to the elementary inequality $\ln(1-\mu)\leq -\mu-\mu^2/2$ by putting
$\mu=1-\lambda$. So, we can finish the proof of Lemma~\ref{L7} by writing
\begin{align*}
\eb\exp\left(\frac\rho4\,\eta_\rho\right)&=\int_{\RR}e^{\,\rho x/4}\;
dF_\rho(x)\\
&=\Bigl[e^{\,\rho x/4}\bigl(F_\rho(x)-1\bigr)\Bigr]_{-\infty}^{+\infty}-\;
\frac\rho4\int_{\RR}e^{\,\rho x/4}\bigl(F_\rho(x)-1\bigr)\;dx\\
&=\frac\rho4\int_{\RR_-}e^{\,\rho x/4}\;dx+
\frac\rho4\int_{\RR_+}e^{\,\rho x/4}\bigl(1-F_\rho(x)\bigr)\;dx\\
&\leq 1+\frac\rho4\int_{\RR_+}e^{-\rho x/4}\;dx=2.
\end{align*}

Now, let us get back to the proof of Lemma~\ref{L5}. Using Lemma~\ref{L7} and
taking into account the stationarity and the independence of the increments of
the process $\ln X_\rho$ on $\RR_+$, we obtain
\begin{align*}
\pb\biggl\{\sup_{y>A} X_\rho(y) > e^{-bA}\biggr\}&\leq e^{\,bA/4}\;
\eb\sup_{y>A} X_\rho^{1/4}(y)\\
&=e^{\,bA/4}\;\eb X_\rho^{1/4}(A)\;
\eb\sup_{y>A}\frac{X_\rho^{1/4}(y)}{X_\rho^{1/4}(A)}\\
&=e^{\,bA/4}\;\eb X_\rho^{1/4}(A)\;\eb\sup_{z>0} X_\rho^{1/4}(z)\\
&=e^{\,bA/4}\;\eb X_\rho^{1/4}(A)\;\eb\sup_{z>0}
\left(\exp\Bigl(\frac\rho4\,\Pi_+(z/\rho)-\frac{z}{4\rho}\Bigr)\right)\\
&=e^{\,bA/4}\;\eb X_\rho^{1/4}(A)\;\eb\exp\left(\sup_{t>0}\Bigl(
\frac\rho4\bigl(\Pi_{\textstyle\frac\rho{e^{\rho}-1}}(t)-t\bigr)\Bigr)\right)\\
&=e^{\,bA/4}\;\eb X_\rho^{1/4}(A)\;\eb\exp\left(\frac\rho4\,\eta_\rho\right)
\leq 2\,e^{\,bA/4}\;\eb X_\rho^{1/4}(A).
\end{align*}

Hence, taking $b\in\left]\,0\,{,}\;3/40\,\right[$, we have
$5b/4\in\left]\,0\,{,}\;3/32\,\right[$ and, using Lemma~\ref{L6}, we finally
get
\begin{align*}
\pb\biggl\{\sup_{y>A} X_\rho(y) > e^{-bA}\biggr\}&\leq
2\,e^{\,bA/4}\,\exp\Bigl(-\frac{5b}{4}A\Bigr)=2\,e^{-bA}
\end{align*}
for all sufficiently small $\rho$ and all $A>0$, and so the first term is
estimated.

The second term can be estimated in the same way, if we show that for all
$\rho>0$ the random variable
$$
\eta_\rho'=\sup_{t\in\RR_+}\bigl(-\Pi_{\lambda'}(t)+t\bigr)
=-\inf_{t\in\RR_+}\bigl(\Pi_{\lambda'}(t)-t\bigr),
$$
where $\Pi_{\lambda'}$ is a Poisson process on $\RR_+$ with intensity
$\lambda'=\rho/(1-e^{-\rho})\in\left]0,1\right[$, verifies
$$
\eb\exp\left(\frac\rho4\,\eta_\rho'\right)\leq 2.
$$

For this, let us remind that according to Pyke~\cite{Pyke} $\bigl($see also
Cram\'er~\cite{Cram}$\bigr)$, $\eta_\rho'$ is an exponential random variable
with parameter $r$, where $r$ is the unique positive solution of the equation
$$
\lambda'(e^{-r}-1)+r=0.
$$
In our case, this equation becomes
$$
\frac{\rho}{1-e^{-\rho}}\,(e^{-r}-1)+r=0,
$$
and $r=\rho$ is clearly its solution. Hence $\eta_\rho'$ is an exponential
random variable with parameter~$\rho$, which yields
$$
\eb\exp\left(\frac\rho4\,\eta_\rho'\right)=\frac43<2,
$$
and so, Lemma~\ref{L5} is proved.

\section{Numerical simulations}
\label{NS}

In this section we present some numerical simulations of the quantities
$B_\rho$, $M_\rho$ and $E_\rho$ for $\rho\in\left]0,\infty\right[$. Besides
giving approximate values of these quantities, the simulation results
illustrate both the asymptotics
$$
B_\rho\sim\frac{B_0}{\rho^2}\,,\quad M_\rho\sim\frac{M_0}{\rho^2}
\quad\text{and}\quad E_\rho\to E_0 \quad\text{as}\quad \rho\to 0,
$$
with $B_0=16\,\zeta(3)\approx19.2329$, $M_0=26$ and
$E_0=8\,\zeta(3)/13\approx 0.7397$, and
$$
B_\rho\to B_\infty,\quad M_\rho\to M_\infty \quad\text{and}\quad E_\rho\to
E_\infty \quad\text{as}\quad \rho\to\infty,
$$
with $B_\infty=1$, $M_\infty=2$ and $E_\infty=0.5$.

First, we simulate the events $x_1,x_2,\ldots$ of the Poisson process $\Pi_+$
$\bigl($with the intensity $1/(e^\rho-1)\bigr)$, and the events
$x'_1,x'_2,\ldots$ of the Poisson process $\Pi_-$ $\bigl($with the intensity
$1/(1-e^{-\rho})\bigr)$.

Then we calculate
\begin{align*}
\zeta_\rho&=\frac{\ds\int_{\RR}x\,Z_\rho(x)\;dx}
{\ds\int_{\RR}\,Z_\rho(x)\;dx}\\
&=\frac {\ds\sum_{i=1}^{\infty}x_i\,e^{\rho i-x_i} + \sum_{i=1}^{\infty}
e^{\rho i-x_i} - \sum_{i=1}^{\infty}x'_i\,e^{\rho-\rho i+x'_i} +
\sum_{i=1}^{\infty}e^{\rho-\rho i+x'_i}}{\ds\sum_{i=1}^{\infty}e^{\rho i-x_i}
+ \sum_{i=1}^{\infty}e^{\rho-\rho i+x'_i}}\\
\intertext{and}
\xi_\rho&=\argsup_{x\in\RR}Z_\rho(x)=\begin{cases}
x_k, &\text{if } \rho k-x_k > \rho-\rho\ell+x'_\ell,\\
-x'_\ell, &\text{otherwise},
\end{cases}
\end{align*}
where
$$
k=\argmax_{i\geq1}\,(\rho i-x_i)\quad\text{and}\quad
\ell=\argmax_{i\geq1}\,(\rho-\rho i+x'_i),
$$
so that
$$
x_k=\argsup_{x\in\RR_+}Z_\rho(x)\quad\text{and}\quad
-x'_\ell=\argsup_{x\in\RR_-}Z_\rho(x).
$$

Finally, repeating these simulations $10^7$ times (for each value of $\rho$),
we approximate $B_\rho=\eb\zeta_\rho^2$ and $M_\rho=\eb\xi_\rho^2$ by the
empirical second moments, and $E_\rho=B_\rho/M_\rho$ by their ratio.

\begin{figure}[!h]
\centering\includegraphics*[width=0.6\textwidth]{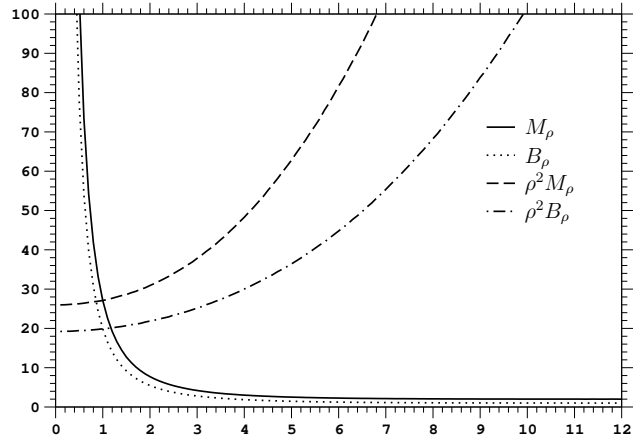}
\caption{$B_\rho$ and $M_\rho$ ($\rho\to 0$ asymptotics)}
\label{fig1}
\end{figure}

\begin{figure}[!h]
\centering\includegraphics*[width=0.6\textwidth]{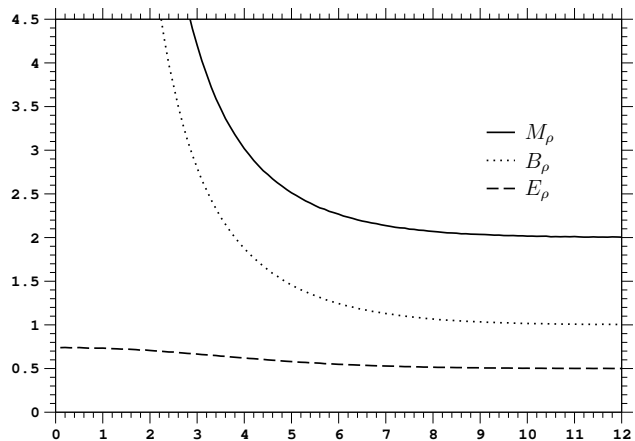}
\caption{\vtop{\hbox{$B_\rho$ and $M_\rho$ ($\rho\to\infty$ asymptotics)}\hbox
{$E_\rho$ (both asymptotics)}}}
\label{fig2}
\end{figure}

The results of the numerical simulations are presented in Figures~\ref{fig1}
and~\ref{fig2}. The $\rho\to 0$ asymptotics of $B_\rho$ and $M_\rho$ can be
observed in Figure~\ref{fig1}, where besides these functions we also plotted
the functions $\rho^2 B_\rho$ and $\rho^2 M_\rho$, making apparent the
constants $B_0\approx19.2329$ and $M_0=26$.

In Figure~\ref{fig2} we use a different scale on the vertical axis to better
illustrate the $\rho\to\infty$ asymptotics of $B_\rho$ and $M_\rho$, as well
as both the asymptotics of $E_\rho$. Note that the function $E_\rho$ appear to
be decreasing, so we can conjecture that bigger is $\rho$, smaller is the
efficiency of the maximum likelihood estimator, and so, this efficiency is
always between $E_\infty=0.5$ and $E_0\approx 0.7397$.

\begin{thebibliography}{99}

\bibitem{CR}\textsc{Chernoff, H} and \textsc{Rubin, H}, ``The estimation of
the location of a discontinuity in density'', \textsl{Proc.\ 3rd Berkeley
Symp.\/}~\textbf{1}, pp. 19--37, 1956.

\bibitem{Cram}\textsc{Cram\'er, H.}, ``On some questions connected with
mathematical risk'', \textsl{Univ.\ California Publ.\ Statist.\/}~\textbf{2},
pp.~99--123, 1954.

\bibitem{DP}\textsc{Deshayes, J.} and \textsc{Picard, D.}, ``Lois
asymptotiques des tests et estimateurs de rupture dans un mod\`ele statistique
classique'', \textsl{Ann.\ Inst.\ H.\ Poincar\'e
Probab.\ Statist.\/}~\textbf{20}, no.~4, pp.~309--327, 1984.

\bibitem{GS}\textsc{Gihman, I.I.} and \textsc{Skorohod, A.V.}, ``\textsl{The
theory of stochastic processes~I.\/}'', Springer-Verlag, New York, 1974.

\bibitem{Gol}\textsc{Golubev, G.K.}, ``Computation of the efficiency of the
maximum-likelihood estimator when observing a discontinuous signal in white
noise'', \textsl{Problems Inform.\ Transmission\/}~\textbf{15}, no.~3,
pp.~61--69, 1979.

\bibitem{HK}\textsc{H\"opfner, R.} and \textsc{Kutoyants, Yu.A.}, ``Estimating
discontinuous periodic signals in a time inhomogeneous diffusion'', preprint,
2009.\\
%
\texttt{http://www.mathematik.uni-mainz.de/\~{}hoepfner/ssp/zeit.html}

\bibitem{IKh2}\textsc{Ibragimov, I.A.} and \textsc{Khasminskii, R.Z.}, ``On
the asymptotic behavior of generalized Bayes' estimator'', \textsl{Dokl.\
Akad.\ Nauk SSSR\/}~\textbf{194}, pp.~257--260, 1970.

\bibitem{IKh4}\textsc{Ibragimov, I.A.} and \textsc{Khasminskii, R.Z.}, ``The
asymptotic behavior of statistical estimates for samples with a discontinuous
density'', \textsl{Mat.\ Sb.\/}~\textbf{87 (129)}, no.~4, pp.~554--558, 1972.

\bibitem{IKh5}\textsc{Ibragimov, I.A.} and \textsc{Khasminskii, R.Z.},
``Estimation of a parameter of a discontinuous signal in a white Gaussian
noise'', \textsl{Problems Inform.\ Transmission\/}~\textbf{11}, no.~3,
pp.~31--43, 1975.

\bibitem{IKh}\textsc{Ibragimov, I.A.} and \textsc{Khasminskii, R.Z.},
``\textsl{Statistical estimation. Asymptotic theory\/}'', Springer-Verlag, New
York, 1981.

\bibitem{KK}\textsc{K\"uchler, U.} and \textsc{Kutoyants, Yu.A.}, ``Delay
estimation for some stationary diffusion-type processes'',
\textsl{Scand.\ J.\ Statist.\/}~\textbf{27}, no.~3, pp.~405--414, 2000.

\bibitem{Kut2}\textsc{Kutoyants, Yu.A.}, ``\textsl{Parameter estimation for
stochastic processes\/}'', Armenian Academy of Sciences, Yerevan, 1980 (in
Russian), translation of revised version, Heldermann-Verlag, Berlin, 1984.

\bibitem{Kut3}\textsc{Kutoyants, Yu.A.}, ``\textsl{Identification of dynamical
systems with small noise\/}'', Mathematics and its Applications~\textbf{300},
Kluwer Academic Publishers Group, Dordrecht, 1994.

\bibitem{Kut}\textsc{Kutoyants, Yu.A.}, ``\textsl{Statistical Inference for
Spatial Poisson Processes\/}'', Lect.\ Notes Statist.~\textbf{134},
Springer-Verlag, New York, 1998.

\bibitem{Kut4}\textsc{Kutoyants, Yu.A.}, ``\textsl{Statistical inference for
ergodic diffusion processes\/}'', Springer Series in Statistics,
Springer-Verlag, London, 2004.

\bibitem{Pfl}\textsc{Pflug, G.Ch.}, ``On an argmax-distribution connected to
the Poisson process'', in \textsl{Proceedings of the Fifth Prague Conference
on Asymptotic Statistics\/}, eds.\ P.~Mandl and H.~Hu\v{s}kov\'a, pp.~123--130,
1993.

\bibitem{Pyke}\textsc{Pyke, R.}, ``The supremum and infimum of the Poisson
process'', \textsl{Ann.\ Math.\ Statist.\/}~\textbf{30}, pp.~568--576, 1959.

\bibitem{RS}\textsc{Rubin, H.} and \textsc{Song, K.-S.}, ``Exact computation
of the asymptotic efficiency of maximum likelihood estimators of a
discontinuous signal in a Gaussian white noise'', \textsl{Ann.\
Statist.\/}~\textbf{23}, no.~3, pp.~732--739, 1995.

\bibitem{ShW}\textsc{Shorack, G.R.} and \textsc{Wellner, J.A.},
``\textsl{Empirical processes with applications to statistics\/}'', John Wiley
\& Sons Inc., New York, 1986.

\bibitem{Ter}\textsc{Terent'yev, A.S.}, ``Probability distribution of a time
location of an absolute maximum at the output of a synchronized filter'',
\textsl{Radioengineering and Electronics\/}~\textbf{13}, no.~4, pp.~652--657,
1968.

\end{thebibliography}
\end{document}